\documentclass{gtart_h}
\def\ifplaintex{\expandafter\ifx\csname documentclass\endcsname\relax}

\def\gtp{{\mathsurround=0pt\it $\cal G\mskip-2mu$eometry \&\ 
$\cal T\!\!$opology $\cal P\!$ublications}}  % GT publications

\def\recd{{\small Received:\qua\receiveddate\ifx\reviseddate\relax
\else\qquad Revised:\qua\reviseddate\fi\par}} 

\def\lognumber#1{\def\thelognumber{#1}}
\def\volumenumber#1{\def\thevolumenumber{#1}}
\def\volumeyear#1{\def\thevolumeyear{#1}}
\def\papernumber#1{\def\thepapernumber{#1}}
\def\pagenumbers#1#2{\def\startpage{#1}\def\finishpage{#2}}
\def\published#1{\def\publishdate{#1}}

\def\received#1{\def\receiveddate{#1}}
\def\revised#1{\def\reviseddate{#1}}
\def\accepted#1{\def\accepteddate{#1}}

\def\asciiauthors#1{\def\theasciiauthors{#1}}
\def\asciiaddress#1{\def\theasciiaddress{#1}}
\def\asciiemail#1{\def\theasciiemail{#1}}

\def\coverauthors#1{\def\thecoverauthors{#1}}
\long\def\asciiabstract#1{\long\def\theasciiabstract{#1}}

\let\\\par\let\thelognumber\relax\let\thevolumenumber\relax
\let\thepapernumber\relax\let\thevolumeyear\relax\let\startpage\relax
\let\finishpage\relax\let\publishdate\relax\let\receiveddate\relax
\let\reviseddate\relax\let\accepteddate\relax\let\theasciititle\relax
\let\theasciiauthors\relax\let\theasciiaddress\relax
\let\theasciiabstract\relax

\let\thecoverauthors\relax\let\theasciiemail\relax

\ifplaintex
\font\logobig=cmssbx10 scaled 3836
\font\logomed=cmssbx10 scaled 2557
\else
\font\logobig=cmssbx10 scaled 4200
\font\logomed=cmssbx10 scaled 2800
\fi

\long\def\makeagttitle{   %%% start of definition of \makeagttitle
\count0=\startpage
\agt\hfill      %   Journal title (top left) 
\hbox to 45truept{\vbox to 0pt{\vglue -13truept{\logomed A\kern -.37em{\logobig 
T}\kern -.38em G}\vss}\hss}
\break
{\small Volume \thevolumenumber\ (\thevolumeyear)
\startpage--\finishpage\nl
Published: \publishdate}

\vglue .25truein

{\parskip=0pt\leftskip 0pt plus
1fil\def\\{\par\smallskip}{\Large\bf\thetitle}\par\medskip} \vglue
0.05truein

{\parskip=0pt\leftskip 0pt plus 1fil\def\\{\par}{\sc\theauthors}
\par\medskip}%
 
\vglue 0.03truein 

{\small\leftskip 25truept\rightskip 25truept{\bf Abstract}\stdspace\theabstract

{\bf AMS Classification}\stdspace\theprimaryclass
\ifx\thesecondaryclass\relax\else; \thesecondaryclass\fi\par
{\bf Keywords}\stdspace \thekeywords\par}\vglue 7truept

}   %%%% end of definition of \makeagttitle

\ifplaintex
\font\phead=cmsl9 scaled 950
\font\pnum=cmbx10 scaled 913
\font\pfoot=cmsl9 scaled 950
\headline{\vbox to 0pt{\vskip -4.5mm\line{\small\phead\ifnum
\count0=\startpage ISSN 1472-2739 (on-line) 1472-2747 (printed)
\hfill {\pnum\folio}\else\ifodd\count0\def\\{ }% 
\ifx\theshorttitle\relax\thetitle\else\theshorttitle\fi\hfill{\pnum\folio}
\else\def\\{ and }{\pnum\folio}\hfill\ifx\theshortauthors\relax\theauthors
\else\theshortauthors\fi\fi\fi}\vss}}
\footline{\vbox to 0pt{\vglue 0mm\line{\small\pfoot\ifnum\count0=\startpage
\copyright\ \gtp\hfill\else
\agt, Volume \thevolumenumber\ (\thevolumeyear)\hfill\fi}\vss}}
\else
\font=cmsl9 scaled 1050
\font\lnum=cmbx10 
\font=cmsl9 scaled 1050
\makeatletter
\def\@oddhead{{\small\lhead\ifnum\count0=\startpage ISSN 1472-2739 
(on-line) 1472-2747 (printed)\hfill {\lnum\number\count0}\else\ifodd\count0
\def\\{ }\ifx\theshorttitle\relax \thetitle \else\theshorttitle\fi\hfill
{\lnum\number\count0}\else\def\\{ and }{\lnum\number\count0}
\hfill\ifx\theshortauthors\relax 
\theauthors\else\theshortauthors\fi\fi\fi}}\def\@evenhead{\@oddhead}
\def\@oddfoot{\small\lfoot\ifnum\count0=\startpage\copyright\ \gtp\hfill\else
\agt, Volume \thevolumenumber\ (\thevolumeyear)\hfill\fi}
\def\@evenfoot{\@oddfoot}
\makeatother
\fi
\let\maketitlepage\makeagttitle
\let\makeshorttitle\maketitlepage
\let\maketitle\maketitlepage

\newwrite\gtoutfile
\long\gdef\makeheadfile{  %%% start of definition of \makeheadfile
{\def\\{, }\def\s{ }
\immediate\openout\gtoutfile head.xxx
\immediate\write\gtoutfile{Proxy-for: \ifx\theasciiauthors\relax
\theauthors\else\theasciiauthors\fi\s<\ifx\theasciiemail\relax\theemail\else\theasciiemail\fi>}
\immediate\write\gtoutfile{\noexpand\\}
\immediate\write\gtoutfile{Authors: \ifx\theasciiauthors\relax
\theauthors\else\theasciiauthors\fi}
{\def\\{ }\immediate\write\gtoutfile{Title: \ifx\theasciititle\relax
\thetitle\else\theasciititle\fi}}
\immediate\write\gtoutfile{Subj-class: GT or SG, GR etc}
\immediate\write\gtoutfile{MSC-class: \theprimaryclass\ifx\thesecondaryclass\relax\else, \thesecondaryclass\fi}
\immediate\write\gtoutfile{Journal-ref: Algebr. Geom. Topol. \thevolumenumber\s
(\thevolumeyear) \startpage-\finishpage}
\immediate\write\gtoutfile{Comments: Published by Algebraic and
Geometric Topology at}
\immediate\write\gtoutfile{\s\s\s  http://www.maths.warwick.ac.uk/agt/AGTVol\thevolumenumber/agt-\thevolumenumber-\thepapernumber.abs.html}
\immediate\write\gtoutfile{\noexpand\\}
\immediate\write\gtoutfile{}
\ifx\theasciiabstract\relax
\immediate\write\gtoutfile{\theabstract}\else
\immediate\write\gtoutfile{\theasciiabstract}\fi
\immediate\write\gtoutfile{}
\immediate\write\gtoutfile{\noexpand\\}
\immediate\write\gtoutfile{}
\immediate\closeout\gtoutfile}}  %%% end of definition of \makeheadfile

\def\maketitlepage{\makeagttitle\makeheadfile}
\let\makeshorttitle\maketitlepage
\let\maketitle\maketitlepage

\lognumber{28}
\volumenumber{5}
\volumeyear{2005}
\papernumber{28}
\pagenumbers{691}{711}
\received{18 October 2004} 
\revised{12 May 2005}
\accepted{27 June 2005}
\published{5 July 2005}

\usepackage{amsmath,amssymb}
\usepackage[arrow,matrix]{xy}

\theoremstyle{plain}
\newtheorem{theorem}[subsection]{Theorem}
\newtheorem{lemma}[subsection]{Lemma}
\newtheorem{prop}[subsection]{Proposition}
\newtheorem{corollary}[subsection]{Corollary}

\theoremstyle{definition}
\newtheorem{remark}[subsection]{Remark}
\newtheorem{definition}[subsection]{Definition}
\newtheorem{example}[subsection]{Example}

\numberwithin{equation}{section}

\newcommand{\A}{{\mathcal A}}
\newcommand{\wA}{\widehat{{\mathcal A}}}
\newcommand{\wB}{\widehat{{\mathcal B}}}

\newcommand{\B}{{\mathcal B}}

\newcommand{\LL}{{\mathcal L}}

\newcommand{\V}{{\mathcal V}}

\newcommand{\Z}{\mathbb{Z}}
\newcommand{\Q}{\mathbb{Q}}

\newcommand{\C}{\mathbb{C}}
\newcommand{\K}{\mathbb{K}}
\newcommand{\PP}{\mathbb{P}}
\newcommand{\FF}{\mathbb{F}}

\newcommand{\tY}{\widetilde{Y}}

\DeclareMathOperator{\Hom}{Hom}
\DeclareMathOperator{\rank}{rank}
\DeclareMathOperator{\im}{im}
\DeclareMathOperator{\coker}{coker}

\DeclareMathOperator{\id}{id}

\DeclareMathOperator{\tor}{Tor}
\DeclareMathOperator{\codim}{codim}
\DeclareMathOperator{\supp}{supp}

\newcommand{\surj}{\twoheadrightarrow}

\begin{document}

\title {Some analogs of Zariski's Theorem\\on nodal line arrangements}

%\shorttitle {Some Analogs of Zariski's Theorem} 

\author{A.D.R. Choudary\\A. Dimca\\\c{S}. Papadima}
\asciiauthors{A.D.R. Choudary, A. Dimca and S. Papadima}
\coverauthors{A.D.R. Choudary\\A. Dimca\\\noexpand\c{S}. Papadima}

\address{Department of Mathematics, Central Washington 
University\\Ellensburg, Washington 98926, USA\qua{\rm and}\\School of Mathematical 
Sciences, GC University Lahore, Pakistan\\\smallskip
\\Laboratoire J.A. Dieudonn\'e, UMR du CNRS 6621\\Universit\'e de
Nice-Sophia-Antipolis, Parc Valrose\\06108 Nice Cedex 02, France\\{\rm and}\\
Inst.\ of Math.\ ``Simion Stoilow", P.O. Box 1-764,
RO-014700 Bucharest, Romania}

\asciiaddress{Department of Mathematics, Central Washington 
University\\Ellensburg, Washington 98926, USA and\\School of Mathematical 
Sciences, GC University Lahore, Pakistan\\Laboratoire 
J.A. Dieudonne, UMR du CNRS 6621\\Universite de
Nice-Sophia-Antipolis, Parc Valrose\\06108 Nice Cedex 02, France\\and\\
Inst. of Math. "Simion Stoilow", P.O. Box 1-764,
RO-014700 Bucharest, Romania}

\gtemail{\mailto{choudary@cwu.edu}, 
\mailto{dimca@math.unice.fr}, 
\mailto{Stefan.Papadima@imar.ro}}

\asciiemail{choudary@cwu.edu, dimca@math.unice.fr, Stefan.Papadima@imar.ro}

\begin{abstract}
For line arrangements in $\PP^2$ with nice combinatorics (in particular, for
those which are nodal away the line at infinity), we prove that the combinatorics
contains the same information as the fundamental group together with the
meridianal basis of the abelianization. We consider higher dimensional analogs
of the above situation. For these analogs, we give purely combinatorial complete
descriptions of the following topological invariants (over an arbitrary field): the
twisted homology of the complement, with arbitrary rank one coefficients; the
homology of the associated  Milnor fiber and Alexander cover, including
monodromy actions; the coinvariants of the first higher non-trivial homotopy
group of the Alexander cover, with the induced monodromy action. 
\end{abstract}

\asciiabstract{%
For line arrangements in P^2 with nice combinatorics (in particular, for
those which are nodal away the line at infinity), we prove that the combinatorics
contains the same information as the fundamental group together with the
meridianal basis of the abelianization. We consider higher dimensional analogs
of the above situation. For these analogs, we give purely combinatorial complete
descriptions of the following topological invariants (over an arbitrary field): the
twisted homology of the complement, with arbitrary rank one coefficients; the
homology of the associated  Milnor fiber and Alexander cover, including
monodromy actions; the coinvariants of the first higher non-trivial homotopy
group of the Alexander cover, with the induced monodromy action.}

 


\maketitle

\section{Introduction}
\label{sec:intro}

Let $\A= \{H_0,...,H_n\}$ be a line arrangement  in the complex projective plane
$\PP^2$
 with complement
$M=M(\A)=\PP^{2}\setminus \cup_{i=0}^n H_i$, and
fundamental group $\pi =\pi_1(M)$. One has the following known result.
\begin{theorem}
\label{Zariski1}

The following statements are equivalent.
\begin{itemize}
\item[\rm(i)] $\A$ has only double points.

\item[\rm(ii)] The fundamental group  $\pi$ is abelian.
\end{itemize}
\end{theorem}

The implication (i) $\implies$ (ii) goes back to Zariski \cite{Z}, Chapter VIII, 
section 2 (see especially his comments on the proof of Theorem 1) and to his attempt
to use it in the proof of his conjecture on the commutativity of the fundamental
group of the
complement of any nodal curve. See Deligne \cite{De} and Fulton \cite{F} for
a proof of this conjecture,
and Dimca \cite[Corollary (4.3.18)]{D1} for details on Zariski's approach.
Concerning this implication, note that Hattori \cite{Hat} proves even more,
namely that the complement $M$ is homotopy equivalent to the 2-skeleton of
the standard $n$-dimensional real torus $T_n$.

The implication (ii) $\implies$ (i) is easy.
Assume that $\A$ has a point $p$ of multiplicity $k>2$ and let $\B$ be the
subarrangement
of $\A$ consisting of all the lines passing through $p$. Then the inclusion
$M(\A) \to M(\B)$ induces an epimorphism at the level of fundamental groups
and $\pi_1( M(\B))$ is a free group $\FF_{k-1}$ on $k-1$ generators.
Indeed, one has $M(\B)=(\C \setminus \{ k-1 ~~ {\rm points}\}) \times \C$,
hence $\pi_1( M(\B))=\FF_{k-1}$.
This implies that $\pi=\pi_1( M(\A)) $ is not abelian, actually, not even nilpotent.

Note also that the implication  (ii) $\implies$ (i) is specific to line 
arrangements.
In fact, there are many non-nodal curves $C \subset \PP^2$ with
$\pi_1(\PP^2 \setminus C)$ commutative, see for such examples \cite[(4.3.8)]{D1}.

In the first part of this paper we investigate to what extent the above Theorem
\ref{Zariski1}
still holds when we consider {\it affine nodal arrangements in the plane} $\C^2$.
In other words, we assume that $H_0$ is the line at infinity $L_{\infty}=\PP^2
\setminus \C^2$
and we ask that $\A$ has only double points in the affine plane $\C^2$. If there is
just one multiple point on the line at infinity (i.e.\ when all the lines in $\A$ 
belong
to the same pencil of lines in $\PP^2$), then $M(\A)$ is of a very special type,
namely a product $(\C \setminus \{ n ~ {\rm points}\}) \times \C$, which is
very easy to treat as an exercise.

In the sequel we assume that there are exactly $r \geq 2$ multiple points on the
line at infinity
$H_0$, and that they are of multiplicities $m_1+1$,...,$m_r+1$. In this situation we
say
that $\A$ is {\it split solvable} of type ${\bf m}=(m_1,m_2,...,m_r)$. 
The arrangements
having this  simple type of combinatorics  belong to the hypersolvable class introduced
by Jambu and Papadima, see \cite[p.1142]{JP}. The type ${\bf m}$ determines
the combinatorics of the arrangement $\A$, i.e.\ its intersection lattice $L(\A)$,
defined as in \cite{OT}. Here is a simple way to think about the type  ${\bf m}$:
there are exactly $r$ directions for the lines in the affine arrangement in $\C^2$
and there are exactly $m_j$ parallel lines having the $j$-th direction as
their common direction, for any $j=1,...,r.$ Requiring only double intersection points,
for the $m_1+...+m_r$ affine lines in $\C^2$, determines then the lattice.

The ideal generalization of  Theorem
\ref{Zariski1}
would be that the following data are equivalent (i.e.\ they determine each other):

\begin{itemize}
\item[(I)] the combinatorics of $\A$ as described by ${\bf m}$;

\item[(II)] the fundamental group $\pi=\pi_1( M(\A)) $.
\end{itemize}

Such a result would be quite interesting, since in general it is known that the
combinatorics
does not determine the fundamental group of a line arrangement, see Rybnikov \cite{Ry}.
However, such a straightforward generalization is not true. Indeed, M. Falk has
produced
in \cite[p.146-147]{Fa} an example of two line arrangements $\A_1$ and $\A_2$
such that $\A_1$ is split solvable of type $(1,2,2)$,
$M(\A_1)$ is homotopy equivalent to $M(\A_2)$, but  $\A_1$ and $\A_2$
are combinatorially distinct.

To state our first result we need some definitions. We assume from now on that all
arrangements are ordered, i.e.\ they come with a numbering of their hyperplanes.
All complex varieties are oriented using the corresponding complex orientations.

The first definition was systematically investigated by Papadima in \cite{P2},
in analogy with classical knot theory.

\begin{definition} \label{defarg}
 Two arrangements $\A= \{H_0,...,H_n\}$ and $\B= \{ {\overline H}_0,...,{\overline
H}_n \}$
in the same ambient projective space $\PP^N$ are said to have
{\it the same oriented topological type} (notation $\A \approx \B$) if
there is an orientation-preserving homeomorphism $f:\PP^N \to
\PP^N$, inducing for all $k=0,...,n$ orientation-preserving homeomorphisms
$f|H_k:H_k \to {\overline H} _k$.

\end{definition}

The second definition we need is the following.

\begin{definition} \label{defmark}
Let $G$ be a group whose abelianization $G_{ab}$ is free abelian of finite rank.
Then a {\it 1-marking} of $G$ is a choice of an ordered basis for $G_{ab}$. If $G$
and $G'$
are two 1-marked groups, we say that $G$ and $G'$ are {\it 1-marked isomorphic} if
there is a group isomorphism $\phi:G \to G'$ preserving in the obvious sense the two
markings.
Notation: $G \approx _1 G'$.

\end{definition}

\begin{example} \label{example1}

If $G=\pi =\pi _1(M(\A))$ is the fundamental group of an ordered projective arrangement
 $\A= \{H_0,...,H_n\}$, then the meridians $\{x_1,...,x_n\}$ give a geometric 1-marking
for $G_{ab}=H_1(M(\A),\Z)$. Here $x_i$ is the homology class of a small loop $\gamma
_i$
going in the positive direction around the hyperplane $H_i$. The fundamental group
of an ordered projective arrangement is always in this paper endowed with this
geometric 1-marking.

Moreover, it is known that  $\A \approx \B$ implies $\pi _1(M(\A)) \approx _1 \pi
_1(M(\B))$
with $\phi = f_{\sharp}$; see  Lemma 2.7 from \cite{P2}.

\end{example}

\begin{remark} \label{rem1}
If $G$ and $G'$ are two finitely presented groups, there is no algorithm
to decide whether they are or not isomorphic. However, if $G$ and $G'$ have
in addition only commutator relations (which is the case for $\pi _1(M(\A))$),
then there is an algorithmic obstruction theory for an analogous 1,2--marked
isomorphism problem, see Papadima \cite[Section 4]{P2}. Moreover, this reduces to
the above $\approx_1$--problem, when $G$ and $G'$ are fundamental groups
of projective arrangements with the same combinatorics, see \cite[Remark 4.7]{P2}.
The second obstruction in this sense plays a key role in Rybnikov's work \cite{Ry}.

\end{remark}

Now we can state our first result. For a type  ${\bf m}=(m_1,m_2,...,m_r)$,
consider the projective space $\PP^r$ with coordinates $z_0,z_1,...,z_r$
and the arrangement $\wA({\bf m})$ therein given by
$$z_0(z_1-z_0)(z_1-2z_0) \cdots (z_1-m_1z_0)(z_2-z_0) \cdots(z_2-m_2z_0) \cdots
(z_r-m_rz_0)=0.$$
Here the hyperplanes in  $\wA({\bf m})$ inherit the order of the corresponding
linear factors in the above product. We set
 $\A({\bf m})=\wA({\bf m}) \cap U$, where $U$ is a general 2-plane in $\PP^r$,
which we identify to $\PP^2$.
It is obvious that the arrangement $\A({\bf m})$ is split solvable of type ${\bf m}$
and that  $M(\wA({\bf m}))$ is isomorphic to the product
$$(\C \setminus \{m_1~ {\rm points}\}) \times \cdots \times (\C \setminus \{m_r~
{\rm points}\}).$$
In particular,  $\wA({\bf m})$ is aspherical. Note also that the combinatorics of
$\A({\bf m})$ is {\it nice}, in the sense of Jiang--Yau \cite[Definition 3.2]{JY2};
see \S~\ref{subsec=pfnice}. 

It turns out that the appropriate generalization of Theorem \ref{Zariski1} actually
holds for
an arbitrary nice combinatorics.

\begin{theorem} \label{thm=nice}
Let $\A$ and $\B$ be line arrangements in $\PP^2$, with $\B$ nice.
The following are equivalent.

\begin{itemize}
\item[\rm(i)]  $L(\A)=L(\B)$ (i.e.\ they have the same combinatorics);

\item[\rm(ii)] $\A \approx \B$ (i.e.\ they have the same oriented topological type);

\item[\rm(iii)] $M(\A) \sim M(\B)$ (i.e.\ their complements are homeomorphic);

\item[\rm(iv)]  $\pi _1(M(\A)) \approx _1 \pi _1(M(\B))$ (i.e.\ their fundamental groups
are $1$--marked isomorphic).
\end{itemize}
\end{theorem}

The implication (iv) $\Rightarrow$ (i), which corresponds to (ii) $\Rightarrow$ (i)
in Zariski's classical set-up from Theorem \ref{Zariski1}, will be proved in the next 
section for 
arbitrary  line arrangements. In fact this will appear as a special case of a result 
on general
hyperplane arrangements, which identifies precisely the amount of combinatorics
determined by the 1-marked
fundamental group of the complement, see Theorem  \ref{newthm}.
The other implications follow from results obtained by
Jiang--Yau in \cite {JY1, JY2}. Theorem \ref{thm=nice} has the following
consequence.

\begin{corollary} \label{mainthm}
Let $\A$ be a line arrangement in $\PP^2$. The following are equivalent.
\begin{itemize}

\item[\rm(i)] $\A$ has the same combinatorics as a split solvable arrangement of type
${\bf m}=(m_1,m_2,...,m_r)$;

\item[\rm(ii)]  $\pi _1(M(\A)) \approx _1 \FF ({\bf m})$, where the product of free groups,
$\FF ({\bf m})=  \FF_{m_1} \times \cdots \times \FF_{m_r}$,
is endowed with the obvious $1$-marking;

\item[\rm(iii)] $\A \approx \A({\bf m})$.
\end{itemize}

\end{corollary}

Following Zariski's ideas, Oka and Sakamoto \cite{OkS} gave a general condition
for the splitting of the fundamental group of the complement of a projective plane
curve
as a direct product. This was used by Fan \cite{Fan} to derive
a sufficient condition (a priori not combinatorial, but verified by the split
solvable class)
under which the fundamental group of a line arrangement is a product of free groups.
However, Fan's paper leaves the converse implication as an open question, see
\cite[p.290]{Fan}.

An important consequence of Zariski's conjecture is that the fundamental group
of a nodal curve has a simple explicit description, involving only the degree list of 
the irreducible components of the curve, see for instance  \cite{D1}.
This viewpoint is pursued in the second part of our paper. More precisely, we extend 
the implication
(i) $\implies$ (ii) from Corollary \ref{mainthm} above to higher dimensional hyperplane 
arrangements with 
simple combinatorics of split solvable type. We show that not only the fundamental group, 
but also a lot of additional important topological information on the complement can be 
{\it explicitly} described in combinatorial terms.
Such computations do not seem at all obvious, for an arbitrary nice combinatorics.

Set $M({\bf m},s):= M(\wA({\bf m}))\cap U^s$, where
$U^s\subset \PP^r$
is a generic $s$--plane, for $2\le s\le r$. Note that 
$M({\bf m},2)= M(\A({\bf m}))$, and
$M({\bf m},r)= M(\wA({\bf m}))= K(\FF ({\bf m}),1)$. We will examine four topological
invariants of $M({\bf m}, s)$, in sections \ref{sec:locsys} and \ref{sec=alex}, 
using tools from \cite{DP1}, \cite{DP2}.

In Section \ref{sec:locsys}, we will look at twisted homology with rank one
coefficients,
and Milnor fiber homology. Both topics are under current intense investigation in
arrangement theory. The important role played in singularity theory by
the homology of the Milnor fiber of a polynomial, together with the monodromy action
on it,
can be traced back to Milnor's seminal book \cite{M}. For the importance of twisted
rank one
homology in arrangement theory, and connexions with hypergeometric integrals, see
the book
by Orlik--Terao \cite{OT2}.

The rank one local systems $\LL_{\rho}$ on $M({\bf m}, s)$ are parametrized by group
homomorphisms,
$\rho \colon \FF({\bf m}) \to \K^*$, where $\K$ is an arbitrary field. The computation
of all Betti numbers of $M({\bf m}, s)$ with coefficients in $\LL_{\rho}$ is carried
out
in Proposition~\ref{prop=rk1}. The result is purely combinatorial, in the sense that it
depends only on ${\bf m}$, $s$, and on the subset of those indices $i$ for which the
restriction of
$\rho$ to $\FF_{m_i}$ is non-trivial.

The {\em Milnor fiber} $F:=f^{-1}(1)$ of a homogeneous degree $d$ polynomial mapping,
$f \colon \C^{s+1} \to \C$, is endowed with a geometric monodromy, induced by
the diagonal action of the cyclic group of complex $d$--roots of unity on $\C^{s+1}$.
Denote by $F({\bf m},s)$ the Milnor fiber of the defining polynomial of the cone
over the
arrangement $\wA ({\bf m})\cap U^s$. Our main result in Section~\ref{sec:locsys} is
Theorem~\ref{thm=milnor}, where we give a complete description of
$H_*(F({\bf m}, s), \K)$, for an arbitrary field $\K$, together with the
$\K\Z_d$--module
structure induced by the geometric monodromy action. Since the result  depends only on
${\bf m}$ and $s$, we deduce in particular that $H_*(F({\bf m}, s), \Z)$ is
torsion--free.
This is to be compared to the examples found by Cohen--Denham--Suciu in \cite{CDS},
where multiarrangements (i.e.\ arrangements defined by non-reduced equations) with
torsion in $H_1(F, \Z)$ are constructed. Note also that the torsion--freeness of
$H_1(F, \Z)$ in the arrangement (i.e.\ reduced) case is still a major open question.

In Section~\ref{sec=alex}, we examine the {\em Alexander cover},
$E({\bf m}, s) \to M({\bf m}, s)$. This is the $\Z$--cover of $M({\bf m}, s)$ with
fundamental group
$N({\bf m}):= \ker ~ (\nu \colon \pi_1(M({\bf m}, s))=\FF({\bf m}) \to \Z)$, where
$\nu$ sends to $1$ all generators of $\FF({\bf m})$. This definition mimicks a
classical
construction from the theory of Alexander invariants of knots and links in $S^3$,
see for
instance Hillman's book \cite{Hill}. In singularity theory, Alexander invariants were
introduced and studied by Libgober, in a series of papers starting with \cite{Li0},
see also
\cite{Li1}, \cite{Li2}, \cite{Li3}, \cite{Li4} as well as  \cite{DN}.

Our first main result in Section~\ref{sec=alex} is Theorem~\ref{thm=alexh}, where we
give
a complete description of the graded $\K\Z$--module $H_*(E({\bf m}, s), \K)$
over an arbitrary field $\K$, only in terms of ${\bf m}$ and $s$.

Our second main result here is related to higher homotopy groups of Alexander covers.
Obviously, $E({\bf m}, r) =K(N({\bf m}), 1)$ is aspherical. Assuming $2\le s < r$,
it turns out that the first higher non-trivial homotopy group of $E({\bf m}, s)$ is
$\pi_s E({\bf m}, s)$. We give a concrete estimate of its non-triviality, in the
following way.
Let $\K$ be an arbitrary field. As is well-known, the $N({\bf m})$--coinvariants,
$$ \big( \pi_s E({\bf m}, s)\otimes \K \big)_{N({\bf m})} \, ,$$
have a natural $\K\Z$--module structure. We describe it completely, in
Theorem~\ref{thm=alexpi},
again solely in terms of ${\bf m}$ and $s$. Finally, we note in Remark~\ref{rk=nonf},
as a result of the aforementioned computation, that $\pi_s E({\bf m}, s)$ is {\em not}
finitely generated over $\Z N({\bf m})$, in spite of the fact that additively
$$ \pi_s E({\bf m}, s) = \pi_s M({\bf m}, s) \, ,$$
and the latter is known to have a finite free resolution over $\Z \FF({\bf m})$,
see~\cite{DP1}.

\begin{remark} \label{rem0}

The statements (i) and (ii) in Theorem \ref{Zariski1} are equivalent to the following.

\begin{itemize}
\item[(iii)] The second Betti number of $M$ is maximal in the set of second Betti numbers
of complements of arrangements of $n+1$ lines in $\PP^2$. Moreover this maximum is
given by
$b_2(M)= {n \choose 2}$.\end{itemize}

The equivalence of (i) and (iii) follows from the description of the
cohomology of the complement of the associated  central arrangement $\A'$ in $\C^3$
given by Orlik and Solomon in \cite{OS}, since the non-nodal points $p$  correspond
exactly
to nontrivial relations in $H^2(M(\A'))$. A further generalization to curve
arrangements
will be  given elsewhere.

\end{remark}

\section{Combinatorics and $1$--marked groups}
\label{sec:proof}

This section is devoted to the proof of  Theorem~\ref{thm=nice} and Corollary~\ref{mainthm}. 
Our key result in this direction is the following.

\begin{theorem} \label{newthm}
Let $\A$ and $\B$ be arbitrary projective arrangements.
Assume that $\pi_1(M(\A)) \approx_1 \pi_1(M(\B))$. Set $p:=\min \{p_{\A}, p_{\B}\}$. 
Then $L_p(\A)= L_p(\B)$. In particular $L_2(\A)= L_2(\B)$.

\end{theorem}

Here $p_{\A}$ stands for the order of $\pi_1$-connectivity  $p(M({\A}))$ of the complement 
$M({\A})$, a homotopy
invariant introduced by Papadima and Suciu in \cite[p.73]{PS}. The notation 
$L_p(\A)= L_p(\B)$ means that $\A$ and $\B$ have the same dependent subarrangements 
of cardinality at most $p+1$.

\begin{example} \label{ex=higher} 

For line arrangements, the above result simply says that the 1-marked group $\pi_1(M({\A}))$ 
determines the intersection lattice $L(\A)$. It seems worth mentioning that this phenomenon 
may also happen in higher dimensions. Indeed, let
$\wA$ and $\wB$ be essential aspherical arrangements in $\PP^s$ and $\PP^t$ respectively. 
Take $(r-1)$--generic proper sections
$\A:= \wA \cap U$ and $\B:= \wB \cap V$, with $r>2$. Then $p_{\A}= p_{\B}=r-1$ by 
\cite[Theorem 18(i)]{DP1}, in particular
$\A$ and $\B$ are non-aspherical. Since $\rank (\A)=\rank (\B)=r$, we deduce 
from Theorem \ref{newthm} that 
$\pi_1(M(\A)) \approx_1 \pi_1(M(\B))$ implies $L(\A)= L(\B)$.

\end{example}

\subsection{Combinatorics and $1$--marked cohomology} \label{subsec=k}

To prove Theorem~\ref{newthm}, we begin by recalling  that the cohomology algebra of 
the complement, together with its natural
$1$--marking, contains the same information as the combinatorics. More precisely, 
we will need the result below, essentially stated
(without proof) by Kawahara in \cite[Lemmas 22 and 23]{Ka}. For the reader's convenience, 
we include a short proof.

\begin{prop} \label{cohoprop}
Let $\A'= \{H_0',...,H_n'\}$ be a central arrangement in $\C^N$. For $i=0,1,...,n$,
let $e_i \in H^1(M(\A'))$ be the cohomology class corresponding to the hyperplane
$H_i'$.
Then the cohomology algebra $H^*(M(\A'))$ and the collection of elements $\{e_0,
...,e_n\}$
determine the combinatorics of $\A'$.
\end{prop}

\begin{proof}
It is enough to notice that, for any  multi-index $i_1 <...<i_p$, one has
$$\codim (H_{i_1}' \cap ... \cap H_{i_p}') <p$$
if and only if
$$ \sum _{j=1,p}(-1)^je_{i_1}\wedge ...\wedge {\hat e_{i_j}} \wedge ... \wedge
e_{i_p} =0$$
in $H^{p-1}(M(\A'))$.

Indeed, the ``only if'' part follows from the description of the cohomology algebra
in \cite{OS}, while the ``if'' part follows by applying Proposition 3.66
in \cite{OT} to the subarrangement of $\A'$ consisting of the hyperplanes
$\{H_{i_1}', ..., H_{i_p}'\}$. 
\end{proof}

\subsection{More on $1$--markings} \label{subsec=more}

As indicated by the previous result, it will be useful to describe the general
relationship
between various types of $1$--markings.

Let $\A \subset \PP^N$ be an arrangement  and let $\A' \subset \C^{N+1}$ be
the associated central arrangement. Then $\pi ' =\pi_1(M(\A'))$ has a natural 1-marking
$\{x_0',...,x_n'\}$, where $x_i' \in H_1((M(\A'))$ is
the homology class of a small loop $\gamma _i'$ going in the positive direction around
the affine hyperplane $H_i'$ corresponding to the projective hyperplane $H_i$.

{}From the natural product decomposition,
$M(\A')=M(\A) \times \C^*$
one gets natural isomorphisms
\begin{equation} \label{ecupi}
\iota (\A) _{\sharp}: \pi_1(M(\A')) {\tilde \to} \pi_1(M(\A)) \times \Z u_0
\end{equation}
and
\begin{equation} \label{ecuho}
\iota (\A) _*: H_1(M(\A')) {\tilde \to} H_1(M(\A)) \oplus \Z u_0
\end{equation}
where $u_0$ is the standard generator of $\pi_1(\C^*)=H_1(\C^*)\simeq \Z$.
The isomorphism $\iota (\A) _*$ sends the classes $x_i'$ to $x_i$ for all $i=1,...,n$
and the sum $x_0'+x_1'+...+x_n'$ to $u_0$, see \cite[p.209]{D2}. It follows that, 
if $\A$ and $\B$ are two projective arrangements  and $\phi: \pi_1(M(\A)) \to \pi_1(M(\B))$
is a marked isomorphism, then
$$\theta = \iota (\B) _{\sharp}^{-1} \circ (\phi \times \id) \circ \iota (\A)
_{\sharp}$$
gives an isomorphism $ \pi_1(M(\A')) \to \pi_1(M(\B'))$ preserving the
geometric 1-markings.

A 1-marking $\{x_1,...,x_n\}$ for a group $G$ induces, for any ring $R$,
an ordered $R$-basis
$\{x_1^*,...,x_n^*\}$  for $H^1(G,R)=\Hom~(H_1(G,\Z),R)$ given by the usual property,
$x_i^*(x_j)=\delta_{ij}$, which will be called a cohomological 1-marking. If
$\phi:G \to G'$ is a 1-marked isomorphism, then $\phi ^*:H^1(G',R) \to H^1(G,R)$
preserves the corresponding cohomological 1-markings.

\subsection{Proof of Theorem~\ref{newthm}}

\noindent Let $\A'$ and $\B'$ be the associated central arrangements. Set 
$\pi '_{\A}:=\pi_1(M(\A'))$,
$K'_{\A}:=K(\pi '_{\A},1)$, denote by $f_{\A}:M(\A') \to K'_{\A}$ the classifying map, 
and likewise for $\B$. Assuming that $\pi_1(M(\A)) \approx_1 \pi_1(M(\B))$, we get 
an algebra isomorphism $H^*(K'_{\A})\simeq H^*(K'_{\B})$, preserving the cohomological 
1-markings, see  subsection \ref{subsec=more}.

We claim that $f_{\A}$ induces isomorphisms $f_{\A}^q:H^q(K'_{\A},\Q) \to H^q(M({\A'}),\Q)$, 
for all $q \leq p(M({\A'}))=p(M({\A}))$, and likewise for $\B$.

To check this, note first that $f_{\A}^q$ is surjective for any $q\geq 0$, since $f_{\A}^1$ 
is by construction an isomorphism and the algebra  $H^*(M({\A'}),\Q)$ is generated by 
$H^1(M({\A'}),\Q)$, \cite{OS}. Next, for $q \leq p_{\A}$, the vector spaces $H^q(K'_{\A},\Q)$ 
and $H^q(M({\A'}),\Q)$ have the same finite dimension, by the very definition of the order 
of $\pi_1$-connectivity $ p_{\A}$, hence our claim. Note that $f_{\A}^1$ takes the 
cohomological 1-marking of
$\pi '_{\A}$ to the $H^1$-basis $\{e_0, ..., e_n\}$ appearing in Proposition \ref{cohoprop}, 
again by construction, and similarly for $f_{\B}^1$. We thus know that the cohomology 
algebras $H^{\leq p}(M({\A'}),\Q)$ and  $H^{\leq p}(M({\B'}),\Q)$ are isomorphic, and 
the distinguished $H^1$-bases are preserved. By Proposition \ref{cohoprop} we get
$L_p(\A)=L_p(\B)$.

Finally, the order of $\pi_1$--connectivity $ p_{\A}$ is at least 2
for any hyperplane arrangement $\A$. This follows from a result of
Randell~\cite{R2}, which says that $f_{\A}^2$ is always an
isomorphism. The proof of Theorem~\ref{newthm} is complete.

\subsection{Proof of Theorem~\ref{thm=nice}} \label{subsec=pfnice}
Jiang and Yau first associate to an arbitrary line arrangement $\A$ a graph
$\Gamma_{\A}$, depending only on $L(\A)$, in the following way. The vertices $\V_{\A}$
are the singular points of $\A$ of multiplicity at least $3$. Two vertices are joined by
an edge if there is a line of $\A$ containing them. For $v\in \V_{\A}$, they define 
a subgraph $\Gamma_{\A}(v)$, having as vertex set, $\V_{\A}(v)$, $v$ and all his
neighbours from $\Gamma_{\A}$. An edge $\{ u, u'\}$ belongs to $\Gamma_{\A}(v)$ if
$u$, $u'$ and $v$ lie on a line of $\A$. The following combinatorial definition appears in
\cite[Definition 3.2]{JY2}.

\begin{definition}
\label{def=nice}

The arrangement $\A$ is nice if there is $\V' \subset \V_{\A}$ such that
$\V_{\A}(u') \cap \V_{\A}(v') =\emptyset$, for all distinct $u', v' \in \V'$,
and with the property that the subgraph $\Gamma' \subset \Gamma_{\A}$, obtained by
deleting the vertex $v'$ and the edges of $\Gamma_{\A}(v')$, for all $v'\in \V'$, is
a forest.
\end{definition}

If $\A$ is split solvable, then all vertices from $\V_{\A}$ lie on the line at
infinity $H_0$. It follows that $\Gamma_{\A}$ is a complete graph, equal to
$\Gamma_{\A}(v)$, for any $v\in \V_{\A}$. In particular $\Gamma'$ is a discrete graph
with no edges, and consequently $\A$ is nice.

We may now give the proofs of Theorem~\ref{thm=nice} and Corollary~\ref{mainthm}.

\medskip 

(i) $\implies$ (ii)\qua If $L(\A)=L(\B)$, then $\A$ and $\B$ are lattice--isotopic; see
\cite[Theorem 3.3]{JY2}. The conclusion follows from Randell's lattice--isotopy theorem
\cite{Ra}.

(ii) $\implies$ (iii)\qua Obvious.

(iii) $\implies$ (i)\qua See \cite{JY1}.

(ii) $\implies$ (iv)\qua This implication is valid for arbitrary arrangements in
$\PP^N$, as
noted in Example~\ref{example1}.

(iv) $\implies$ (i)\qua See our Theorem~\ref{newthm}.

\subsection{Proof of Corollary~\ref{mainthm}} \label{subsec=pfcorss}
To get Corollary~\ref{mainthm} from Theorem~\ref{thm=nice}, all we need is to
prove the marked isomorphism
$\pi_1(M(\A({\bf m}))) \approx_1 \FF({\bf m})$. This in turn follows from  our
definition
of  $\A({\bf m})$ and from Zariski's Theorem on hyperplane sections, see for instance
\cite[p.26]{D1}, saying that the inclusion $M(\A({\bf m}) ) \to M(\wA({\bf m}) )$
is a 2-equivalence. Therefore
$$\pi_1(M(\A({\bf m}) ))  \approx_1  \pi_1(M(\wA({\bf m}) ))  \approx_1
\FF_{m_1} \times \cdots \times \FF_{m_r}.$$

\section{Combinatorial formulae for twisted homology}
\label{sec:locsys}

We know from Corollary~\ref{mainthm} that the complement of a split--solvable line
arr\-angement
of type ${\bf m}$, $M(\A)$, is homeomorphic to the combinatorial model, $M(\wA({\bf
m}))\cap U^2$,
where $U^2\subset \PP^r$ is a generic 2--plane. Our next results will describe
various topological
invariants of $M(\A)$, in terms of ${\bf m}$.

\subsection{The set-up} \label{subsec=frame}

It turns out that this may be done, at no extra cost, for all generic $s$--plane 
sections,
$M({\bf m},s):= M(\wA({\bf m}))\cap U^s$, with $2\le s\le r$. Set $M:= M({\bf m},s)$
and
$\pi:= \pi_1(M)$. Theorem 18 from \cite{DP1} guarantees the existence, up to homotopy,
of a minimal CW--structure $Y$ on $M(\wA({\bf m}))=K(\FF({\bf m}), 1)$, whose
$s$--skeleton
is homotopy equivalent to $M$: $M \cong Y^{(s)}$.

We aim at computing twisted homology, $H_*(M, S)$, where $S$ is a (right)
$\Z\pi$--module.
For this, we may use the $\pi$--equivariant chain complex of the universal cover of
$Y$,
\begin{equation} \label{eq=chainy}
C_{\bullet}(\tY) \, \colon \, 0\to \Z\pi \otimes C_r \to \cdots \to \Z\pi \otimes C_s
\stackrel{d_s}{\rightarrow} \Z\pi \otimes C_{s-1} \to \cdots \, ,
\end{equation}
where $C_k$, $0\le k \le r$, denotes the free abelian group generated by
the $k$--cells of $Y$. Indeed, one knows (see e.g. \cite[Chapter VI]{W}) that
\begin{equation} \label{eq=wh}
H_*(M,S) = H_*(S \otimes_{\Z\pi} C_{\le s}(\tY))\, .
\end{equation}
Since $s \ge 2$, $\pi =\pi_1(M(\wA({\bf m}))$, which is a product of free groups,
$\FF({\bf m}) =\FF_{m_1} \times \cdots \times \FF_{m_r}$. Consider now the well-known
standard minimal resolution of $\Z$ over $\Z\FF(x_1, \dots ,x_m)$, 
where $\FF(x_1, \dots,x_m)$ is the free group on the letters $x_1, \dots ,x_m$, 
see for instance
\cite[pp.196-197]{HiS}. Use then tensor products of resolutions, as in
\cite[p.222]{HiS},
to obtain a free resolution of $\Z$ over $\Z\pi$,
\begin{equation} \label{eq=chainpi}
C_{\bullet}(\pi) \, \colon \, 0\to \Z\pi \otimes C_r \to \cdots \to \Z\pi \otimes C_s
\stackrel{\partial_s}{\rightarrow} \Z\pi \otimes C_{s-1} \to \cdots
\end{equation}
Note that both \eqref{eq=chainy} and \eqref{eq=chainpi} are minimal free resolutions
of $\Z$ over $\Z\pi$. Consequently, $C_k =\tor_k^{\Z\pi} (\Z, \Z)$, for all $k$. In
more concrete terms, we see from the construction of resolution \eqref{eq=chainpi} that
the graded abelian group $C_{\bullet}:= \oplus_{k=0}^r C_k$ equals
\begin{equation} \label{eq=tensorad}
C_{\bullet}= \bigotimes_{i=1}^r C_{\bullet}^i \, .
\end{equation}
Here $C_{\bullet}^i$ is concentrated in degrees $\bullet =0$ and $1$, for all $i$,
with graded
pieces $\Z$ and $\Z^{m_i}$ respectively. Denoting by $\{ v_1^i, \dots, v_{m_i}^i \}$
the
standard basis of $\Z^{m_i}$, and setting $\FF_{m_i}= \FF(x_1^i, \dots, x_{m_i}^i)$,
we may completely describe the resolution \eqref{eq=chainpi}, as follows:
\begin{equation} \label{eq=expl}
\partial_k (1\otimes v_{j_1}^{i_1}\otimes \cdots \otimes v_{j_k}^{i_k}) =
\sum_{p=1}^k (-1)^{p-1} (x_{j_p}^{i_p} -1)\otimes   v_{j_1}^{i_1}\otimes \cdots
\widehat{v_{j_p}^{i_p}} \cdots \otimes v_{j_k}^{i_k}\, .
\end{equation}
We may now explain our computational strategy. For $S$ a principal ideal domain, with
$\Z\pi$--module structure given by a change of rings, $\Z\pi \to S$, we may use
\cite[Lemma 2.10]{DP2} to infer that the $S$--chain complexes $S\otimes_{\Z\pi}
C_{\bullet}(\tY)$
and $S\otimes_{\Z\pi} C_{\bullet}(\pi)$ are actually isomorphic,
not just chain--homotopy equivalent. This will enable us to replace $C_{\le s}(\tY)$
in equation \eqref{eq=wh} by a simpler, explicit object, namely
$C_{\le s}(\pi)$ from \eqref{eq=chainpi}.

One more remark will be useful. Each graded abelian group $C_{\bullet}^i$ appearing
in \eqref{eq=tensorad} is actually a chain complex, with differential
$\partial' \colon \Z^{m_i} \to \Z$ given by
\begin{equation} \label{eq=dif}
\partial' (v_j^i)= 1\, , \quad {\rm for} \quad j=1, \dots, m_i \, .
\end{equation}
We may thus view $C_{\bullet}$ as a tensor product chain complex,
\begin{equation} \label{eq=tensordif}
(C_{\bullet}, \partial') = \bigotimes_{i=1}^r (C_{\bullet}^i, \partial') \, .
\end{equation}

\subsection{Local systems} \label{subsec=loc}

A rank one local system $\LL_{\rho}$ of $\K$-vector spaces on
$M=M({\bf m}, s)$
corresponds to a representation $\rho:\FF({\bf m}) \to \K^*$, where
$\rho=(\rho_1,...,\rho_r)$
with $\rho_i: \FF_{m_i} \to \K^*$  a representation of the free group
$ \FF_{m_i}$.
We will apply Theorem 4.5 (2) in \cite{DP2} to express the dimensions
of the cohomology groups $H^*(M({\bf m}, s),\LL_{\rho})$ in terms of
the type ${\bf m}$, $s$  and the {\it support of the representation}
$$\supp~(\rho)=\{1 \leq i \leq r ~|~ \rho _i \ne 1\}\, .$$
This extends the results of Hattori from \cite{Hat}, where
the case $m_1=...=m_r=1$ is treated.

For the above representation $\rho$  we define its Poincar\'e
polynomial by the formula
\begin{equation} \label{Ppoly}
P_{\rho}(t)= \sum_{q \geq 0}b_q(\rho)t^q
\end{equation}
where
\begin{equation} \label{Betti}
b_q(\rho)=\dim_{\K}\tor_q^{\Z\FF({\bf m})}(\Z,\K_{\rho})\, .
\end{equation}
To apply the aforementioned theorem, we will need the following standard computation.

\begin{lemma} \label{Betticor}
Let $\rho:\FF({\bf m}) \to \K^*$ be an arbitrary representation. Then
$$P_{\rho}(t)=\prod_{i \notin \supp(\rho)}(1+m_it) \prod_{i \in
\supp(\rho)}(m_i-1)t\, .$$

\end{lemma}

\begin{proof}
Firstly, one has the following equality.
 $$P_{\rho}(t)=\prod_{i=1,r}P_{\rho_i}(t)\, .$$
This equality follows via the K\"unneth formula from resolution
\eqref{eq=chainpi} or, alternatively, from the K\"unneth formula
for cohomology with sheaf coefficients, see \cite[p.117]{D2}.
The case of a single free group is immediate, using the standard minimal resolution.
\end{proof}

Denoting by $\sigma_i$, $0\le i\le r$, the value of the $i$-th elementary symmetric
function on ${\bf m}= (m_1,\dots, m_r)$, we may now state our next result.

\begin{prop} \label{prop=rk1}
Lat $M({\bf m}, s)$ denote the $s$--generic section of $M(\wA({\bf m}))$,
for $2\le s\le r$.
Let $\rho \colon \pi_1(M({\bf m}, s))=\FF({\bf m}) \to \K^*$ be an arbitrary
rank one local system  of $\K$--vector spaces on $M({\bf m}, s)$. Then
\[
\dim_{\K} H_j(M({\bf m}, s), \LL_{\rho})=
\left\{
\begin{array}{cll}
b_j(\rho) \, , & \text{for} & j<s\, ;\\
\sigma_s + \sum_{i=0}^{s-1} (-1)^{i+s} [\sigma_i -b_i(\rho)]\, ,
& \text{for} & j=s\, ;\\
0\, , & \text{for} & j>s\, ,
\end{array}
\right.
\]
where the twisted Betti numbers $\{ b_i(\rho) \}_{0\le i\le r }$ are computed in
Lemma~\ref{Betticor}.
In particular, for all $0 \leq j \leq s$, one has
$$\dim_{\K} H_j(M({\bf m}, s), \K)=\sigma_j.$$

\end{prop}

\begin{proof}
As recalled in \S\ref{subsec=frame}, $M(\wA({\bf m}))\cong Y$ and
$M({\bf m},s)\cong Y^{(s)}$, where $Y$ is a minimal complex. Obviously, the
Poincar\' e polynomial of $\wA({\bf m})$ equals $\prod_{i=1,r} (1+m_i t)$.
Therefore, $\chi (M({\bf m},s))=\sum_{j=0,s} (-1)^j \sigma_j$. Everything
follows now from \cite[Theorem 4.5 (2)]{DP2}.
\end{proof}

In the case of rank one twisted homology, where $S=\K$, things are easy, since
we only need to perform computations involving the fundamental group, and we may
finish by using an Euler characteristic argument. In the sequel, $S$ will be the
$\K$--algebra of a cyclic group. To identify the $S$--module structure of
$$H_*(M,S)=H_*(S\otimes_{\Z\pi} C_{\le s}(\tY))$$
(see \cite[Ch.VI]{W}), we need more, namely to replace
$S\otimes_{\Z\pi} C_{\bullet}(\tY)$ by $S\otimes_{\Z\pi} C_{\bullet}(\pi)$, as
explained
in \S\ref{subsec=frame}.

\subsection{The homology of the associated Milnor fiber}
\label{subsec=milnorh}

Let $\A'$ be a central arrangement of $n+1$ hyperplanes in $\C^{N+1}$. Let $f=0$
be a reduced defining equation for $\A'$ in $\C^{N+1}$. Then $f$ is a homogeneous
polynomial of degree $d=n+1$. The affine smooth hypersurface $F$ given by the
equation $f=1$ in
$\C^{N+1}$ is called the associated {\em Milnor fiber} of the  central arrangement
$\A'$. This space comes equipped with an order $d$ diffeomorphism $h$, the monodromy,
given by $h(x)=\lambda \cdot x$, where $\lambda=\exp(2 \pi \sqrt{-1}/d)$. This
diffeomorphism gives rise to an important (left) $R\Z_d$--module structure on
$H_*(F,R)$,
for any commutative ring $R$. Our next goal is to describe completely this
structure, for
the generic section $\A =\wA({\bf m})\cap U^s$ and an arbitrary field $R=\K$, in terms
of ${\bf m}$ and $s$.

Keeping the previous notation, we may do this as follows. Denote by $\nu\colon \pi
\to \Z$
the (unique) abelian character which takes the value $1$ on the distinguished basis of
$\pi_{ab}=\FF({\bf m})_{ab}$. Note that $\K\Z_d$ is a quotient ring of the principal
ring $\K\Z$,
with $\Z\pi$--module structure coming from the change of rings
$ \Z\pi \stackrel{\nu}{\to} \K\Z \surj \K\Z_d$.
One knows (see \cite{CS}  and \cite[(2.3.4)]{D2}) that $H_*(F, \K)=H_*(M,
\K\Z_d)$,
as $\K\Z_d$--modules. It follows from \S\ref{subsec=frame} that
\begin{equation} \label{eq=change}
H_*(M, \K\Z_d)= H_*(\K\Z_d\otimes_{\Z\pi} C_{\le s}(\pi))\, ,
\end{equation}
as $\K\Z_d$--modules.

We now remark that the differential of the chain complex $\K\Z_d
\otimes_{\Z\pi}C_{\bullet}(\pi)$ equals
\begin{equation}
\label{eq=key1}
\{ (\tau -1)\otimes \partial'_k \colon \K\Z_d \otimes C_k \rightarrow \K\Z_d \otimes
C_{k-1} \}_k\, ,
\end{equation}
where $\partial'$ is described in \eqref{eq=tensordif} and $\tau -1$ is the
multiplication operator
by $(\tau -1)\in \K\Z_d$, $\tau$ being the canonical generator of $\Z_d$. This follows
from \eqref{eq=expl}, together with \eqref{eq=dif}.

It follows from \eqref{eq=key1} that
\begin{equation}
\label{eq=z}
\ker~\big( (\tau -1)\otimes \partial'_k \big) = \big( \K\Z_d \otimes
\ker~(\partial'_k \otimes \K)\big)
\bigoplus \big( (\K\Z_d)^{\Z} \otimes \im~(\partial'_k \otimes \K)\big)\, ,
\end{equation}
where $(\K\Z_d)^{\Z}:= \ker~(\tau -1)=\K$ with trivial $\K\Z_d$--action, and
\begin{equation}
\label{eq=b}
\im~\big( (\tau -1)\otimes \partial'_{k+1} \big) = (\tau -1) \K\Z_d \otimes
\im~(\partial'_{k+1} \otimes \K)\, .
\end{equation}
We are thus led to compute the homology of $(C_{\bullet}, \partial')$.

\begin{lemma}
\label{lem=chom}
For any coefficient ring $R$, $H_*(C_{\bullet}, R)$ is concentrated in degree $*=r$,
where it is free of rank $(m_1-1)\cdots(m_r-1)$.

\end{lemma}

\begin{proof}
Follows from equations \eqref{eq=dif}--\eqref{eq=tensordif}, via the K\" unneth
formula.
\end{proof}

To state our main result in this section, we set
\begin{equation}
\label{eq=defz}
z_j~ :=~ \sigma_j +\sum_{i=o}^{j-1}(-1)^{i+j} \sigma_i\, , \quad {\rm for} \quad
0\le j\le r \,.
\end{equation}

\begin{theorem}
\label{thm=milnor}
Let $F({\bf m}, s)$ denote the Milnor fiber of the generic $s$--section,
$\wA({\bf m})\cap U^s$. Set $d=1+\sigma_1$. Then
\[
H_j(F({\bf m}, s), \K)=
\left\{
\begin{array}{cll}
\K^{z_j +z_{j-1}} \, , & \text{for} & j<s\, ;\\
\K^{z_{s-1}} \oplus (\K\Z_d)^{z_s} \, ,
& \text{for} & j=s\, ;\\
0\, , & \text{for} & j>s\, ,
\end{array}
\right.
\]
as $\K\Z_d$--modules, where $\K$ is an arbitrary field. In particular,
$H_*(F({\bf m}, s), \Z)$ is torsion--free and the integral monodromy action is 
trivial in degrees $<s$.
\end{theorem}

\begin{proof}
An Euler characteristic argument may be used to infer from Lemma~\ref{lem=chom}
that $\dim_{\K} \ker~(\partial'_j\otimes \K) =z_j$, for $j=0,\dots,r$. Recall that
the $\K\Z_d$--module $H_*(F({\bf m}, s), \K)$ may be computed from \eqref{eq=change}.
The first assertion of the theorem follows then from equations
\eqref{eq=key1}--\eqref{eq=b},
together with Lemma~\ref{lem=chom}. Since the $\K$--Betti numbers of $F({\bf m}, s)$
are
independent of $\K$, the $\Z$--homology has no torsion.
\end{proof}

\begin{remark}
\label{rk=zpos}
From \eqref{eq=defz} we see that $z_r=(m_1-1)\cdots (m_r-1)$ is positive, unless
there is
some $m_i$ equal to $1$. We note that always $z_j>0$, if $j<r$. This may be seen
for instance as follows. If $z_j=0$, then $\partial'_{j+1}=0$. On the other hand,
plainly
$C_i\neq 0$, for any $i\le r$. Pick then any basis element $v\in C_{j+1}$ and compute
$\partial'_{j+1}(v)$ using \eqref{eq=dif}--\eqref{eq=tensordif} to obtain a
contradiction.
\end{remark}

\begin{remark}
\label{rk=mtame}
The triviality of the monodromy action on $H^1(F({\bf m}, 2), \C)$ was
obtained by the second author along different lines, using relations to the topology of a polynomial mapping
$\C^2 \to \C$, see \cite[Corollary 3.3]{D3}. This was the starting point of our
interest
in the class of split solvable line arrangements, and their higher dimensional analogs.

The same triviality of the complex monodromy action holds in the dual case when a line 
arrangement $\A$ has a line along which the only singularities are nodes, see \cite{Li4} 
and \cite{CDO} for higher dimensional analogs of this result.

Note also that, for Milnor fibers associated to projective hypersurfaces, the triviality 
of the complex monodromy may occur even when the integral monodromy is not trivial, 
see \cite{Sie}.
\end{remark}

\section{Homology and homotopy groups of infinite cyclic covers}
\label{sec=alex}

Let $M$ be a connected complex. Any $1$--marking of $\pi:=\pi_1(M)$ gives rise to an
abelian
character, $\nu \colon \pi \to \Z$, which takes the value $1$ on all elements of the
distinguished basis of $\pi_{ab}$. The $\Z$--cover associated to $\nu$, $q\colon E
\to M$,
will be called the {\em Alexander cover}, by analogy with the classical Alexander
polynomial
theory (in one variable) for links. In this section, we will take $M=M({\bf m}, s)$,
like in
Section~\ref{sec:locsys}, endowed with the canonical $1$--marking of $\pi =\FF({\bf
m})$, and
compute two basic topological invariants of $E$.

\subsection{Homology of the Alexander cover}
\label{subsec=halex}

The tools we have developed so far enable us to give the following complete
description of
the homology with field coefficients, which involves only ${\bf m}$ and $s$, as before.

\begin{theorem} \label{thm=alexh}
Denote by $E({\bf m}, s)$, $2\le s\le r$, the Alexander cover of the complement of
the generic
$s$--section of $\wA({\bf m})$. Then
\[
H_j(E({\bf m}, s), \K)=
\left\{
\begin{array}{cll}
\K^{z_j} \, , & \text{for} & j<s\, ;\\
(\K\Z)^{z_s} \, ,
& \text{for} & j=s\, ;\\
0\, , & \text{for} & j>s\, ,
\end{array}
\right.
\]
as $\K\Z$--modules (where $\{ z_j \}_{0\le j\le r}$  are defined in
\eqref{eq=defz}), for
any field $\K$.

\end{theorem}

\begin{proof}
It follows from~\cite[Ch.VI]{W} that
\begin{equation}
\label{eq=e1}
H_*(E, \K)=H_*(\K\Z\otimes_{\Z\pi} C_{\le s}(\tY))\, ,
\end{equation}
as $\K\Z$--modules, where $C_{\bullet}(\tY)$ is as in \eqref{eq=chainy} and the
$\Z\pi$--module structure on $\K\Z$ is given by $\nu\colon \Z\pi \to \K\Z$. As
explained in
\S\ref{subsec=frame}, equation \eqref{eq=e1} above may be replaced by the simpler
explicit form
\begin{equation}
\label{eq=e2}
H_*(E, \K)=H_*(\K\Z\otimes_{\Z\pi} C_{\le s}(\pi))\, ,\quad {\rm over} \quad \K\Z \, .
\end{equation}
From now on, the computation is almost identical to the one from the proof of
Theorem~\ref{thm=milnor} from \S\ref{subsec=milnorh}, with $\K\Z_d$ replaced by $\K\Z$,
except the fact that in the case of the Alexander cover we have
$(\K\Z)^{\Z}=\ker~(\tau -1)=0$.
\end{proof}

\subsection{Homotopy of the Alexander cover}
\label{subsec=pialex}

Let us analyze now the homotopy groups of the Alexander cover,
$q\colon E({\bf m}, s)\to M({\bf m}, s)$. Set
$N({\bf m}):= \pi_1(E({\bf m}, s))$. Clearly,
$N({\bf m})=\ker~(\nu \colon \FF({\bf m}) \surj \Z)$. Since
$M({\bf m}, r)=M(\wA({\bf m}))$ is aspherical, $E({\bf m}, r)=K(N({\bf m}), 1)$. We
will
thus assume in the sequel that $2\le s\le r-1$, and suppress ${\bf m}$ and $s$ from
notation.

It follows from \cite[Theorem 18$(ii)$]{DP1} that $\pi_sE=\pi_sM$ is the first
higher non-trivial
homotopy group of $E$. Moreover, $\pi_sM$ has a finite free, minimal,
$\Z\pi$--resolution,
obtainable from \eqref{eq=chainy}. In particular, $\pi_s M$ has  the following
finite, minimal,
$\Z\pi$--presentation
\begin{equation}
\label{eq=mpres}
\pi_s M=\coker~\big( d_{s+2} \colon \Z\pi \otimes C_{s+2} \rightarrow \Z\pi \otimes
C_{s+1} \big)\, .
\end{equation}
Taking $\pi$--coinvariants in equation~\eqref{eq=mpres} (that is, passing from
the group $\pi_s M$ to its coinvariants
$(\pi_s M)_{\pi}:= \pi_s M/I\pi \cdot \pi_s M$, where $I\pi \subset \Z\pi$ is the
augmentation ideal),
we find out from minimality that
\begin{equation}
\label{eq=mco}
(\pi_s M)_{\pi}= C_{s+1}\, ,
\end{equation}
a finitely generated free abelian group, with combinatorially determined rank.

Guided by this result, we are going to look in the sequel at $\pi_sE$ as $\Z
N$--module,
more precisely at the $N$--coinvariants,
\[
\big( \pi_s E\otimes \K \big)_N := (\pi_sE \otimes \K)/(IN \otimes \K)\cdot (\pi_sE
\otimes \K)\, ,
\]
with arbitrary field coefficients. Note that $N$ is a normal subgroup of $\pi$, and
the $\Z N$--module
structure on $\pi_s E$ comes from the $\Z\pi$--module structure of $\pi_s M$, by
restriction of scalars.
Consequently, $(\pi_s E\otimes \K)_N$ has a naturally induced $\K \Z$--module
structure.

We will completely describe this structure, again solely in terms of ${\bf m}$ and
$s$. Here is our second
main result in this section.

\begin{theorem}\label{thm=alexpi}
Let $E({\bf m}, s)$, $2\le s <r$, be the total space of the Alexander cover of
$M({\bf m}, s)$, with
respect to the canonical $1$--marking of $\pi_1(M({\bf m}, s))=\FF({\bf m})$.
Set $N({\bf m})=\pi_1(E({\bf m}, s))$. Let $\K$ be an arbitrary field. Then
the $\K \Z$--module structure of the $N({\bf m})$--coinvariants of the first higher
non-trivial
homotopy group of $E({\bf m}, s)$ is given by
\[
\big( \pi_s E({\bf m}, s)\otimes \K \big)_ {N({\bf m})} =
\left\{
\begin{array}{cll}
(\K\Z)^{z_s} \oplus \K^{z_{s+1}} \, , & \text{for} & s<r-1\, ;\\
(\K\Z)^{\sigma_r} \, ,
& \text{for} & s=r-1\, ,
\end{array}
\right.
\]
where $\{ \sigma_j:= \sigma_j({\bf m}) \}_{0\le j \le r}$ are the elementary symmetric
functions evaluated at ${\bf m}$, and $\{ z_j \}_{0\le j \le r}$ are defined by
\eqref{eq=defz}.
\end{theorem}

\begin{proof}
Firstly, we infer from \eqref{eq=mpres} that
\begin{equation}
\label{eq=a1}
\big( \pi_s E\otimes \K \big)_N =\coker~ \big( \K\Z \otimes_{\Z\pi} d_{s+2} \colon
\K\Z \otimes C_{s+2} \rightarrow \K\Z \otimes C_{s+1} \big) \, ,
\end{equation}
as $\K\Z$--modules.

Once again, we may replace $\K\Z \otimes_{\Z\pi} C_{\bullet}(\tY)$ by
$\K\Z \otimes_{\Z\pi} C_{\bullet}(\pi)$, as explained in \S\ref{subsec=frame} and
\S\ref{subsec=milnorh}, to arrive at
\begin{equation}
\label{eq=a2}
\big( \pi_s E\otimes \K \big)_N =\coker~ \big( (\tau -1)\otimes \partial'_{s+2} \colon
\K\Z \otimes C_{s+2} \rightarrow \K\Z \otimes C_{s+1} \big) \, ,
\end{equation}
over $\K\Z$.

All assertions to be proved follow now from \eqref{eq=a2} above, by using
Lemma~\ref{lem=chom}
and the $\K\Z$ analog of \eqref{eq=b}.
\end{proof}

\begin{remark} \label{rk=nonf}
Upon applying Remark~\ref{rk=zpos} to the above theorem, we see that $\pi_sE({\bf m}, s)$ is
{\em not} finitely generated over the group ring $\Z\pi_1(E({\bf m}, s))$, if $s<r$. This is
in striking contrast with the strong finiteness property of $\pi_sM({\bf m},
s)=\pi_sE({\bf m}, s)$
over $\Z\pi_1(M({\bf m}, s))$, recalled at the beginning of \S\ref{subsec=pialex}. The
non--finiteness of $\pi_{r-1}E({\bf m}, r-1)$ over $\Z N({\bf m})$, in the case when
$m_i>1$, for all $i$, is related to a certain non--finiteness homological property  of
the Bestvina--Brady groups $N({\bf m})$, see \cite{BB}, as first noticed by
Stallings \cite{St}.
\end{remark}

\medskip
{\bf Acknowledgements}\qua We are grateful to Alex Suciu for numerous, very
useful
discussions. This work started while
the third author was visiting the University of Bordeaux, and it was finished
during his stay at the  Mathematical Sciences Research Institute, Berkeley.
He is thankful to both institutions, for
support and excellent working facilities.

A.D.R. Choudary and A. Dimca have been partially funded for this research by a grant from Higher Education
Commission, Pakistan.
S. Papadima has been partially supported by CERES grant
4-147/12.11.2004 of the Romanian Ministry of
Education and Research.

\Addresses\recd

\end{document}